\documentclass[twoside]{article}
\usepackage{upref}
\usepackage{amssymb,amsfonts, amsmath}
\usepackage[all]{xy}
\begin{document}
\title{\vspace{-1in}\parbox{\linewidth}{\footnotesize\noindent}
The fixed point property for a class of nonexpansive maps in  $L\sp{\infty}(\Omega,\Sigma,\mu)$\\
\vspace{.5cm}
\large Cleon S. Barroso
\vspace{-.5cm}
\author{\footnotesize{\it Departamento de Matem\'atica, 
Universidade Federal do Cear\'a},  \\
\footnotesize\it Campus do Pici, \sl B.914, 60455-760, \sl Fortaleza, CE, Brazil
}
\thanks{ {\em E-mail address:} cleonbar@mat.ufc.br
\hfill\break\indent
{\em Key words:} Fixed point property, weakly compact convex, $L\sp{\infty}(\Omega,\Sigma,\mu)$.
}}
\date{}
%
\maketitle 

\newtheorem{theorem}{Theorem}[section]
\newtheorem{remark}{Remark}
\newtheorem{proposition}{Proposition}[section]
\newtheorem{example}{Example}
\newtheorem{corollary}{Corollary}[section]
\numberwithin{equation}{section}

\allowdisplaybreaks

\begin{abstract}
For a finite and positive measure space $(\Omega,\Sigma,\mu)$ and any weakly compact convex subset of $L\sp{\infty}(\Omega,\Sigma,\mu)$ a fixed point theorem for a class of nonexpansive self-mappings is proved. An analogous result is obtained for the space $C(\Omega)$. An illustrative example is given.
\end{abstract}

\section{Introduction}\label{sec:intro}

The problem if every weakly compact convex $K$ of a Banach space $X$ has the fixed point property for nonexpansive mappings has been intensivelly studied for several years, see \cite{1,2,3,4}, and the references therein. Let $(\Omega,\Sigma,\mu)$ be a finite and positive measure space. It is well-known that closed balls and weak$\sp{\ast}$ compact subsets of $L\sp{\infty}(\Omega,\Sigma,\mu)$ always have the fixed point property with respect to nonexpansive self-mappings, see \cite{5}. However, even in particular cases no result seems to be known for weakly compact convex subsets of $L\sp{\infty}(\Omega,\Sigma,\mu)$.
\par
In this paper we study this problem for an special class of nonexpansive maps in $L\sp{\infty}(\Omega,\Sigma,\mu)$, namely, the class of strongly nonexpansive maps (see definition below). We will show that any weakly compact convex subset of $L\sp{\infty}(\Omega,\Sigma,\mu)$ have the fixed point property for strongly nonexpansive maps. Although our result is restricted to a class of nonexpansive mappings, it indicates which type of maps must be rejected in case that intends to show that $L\sp{\infty}(\Omega,\Sigma,\mu)$ fails the fixed point property for such sets. 
Moreover, as a consequence of the ideas developed here a similar result can be proved for the space $C(\Omega)$, where $\Omega$ is a compact Hausdorff space.
As a simpler application, we use this last result to shows that the assumption of weak compactness in the theorem of Arino, Gautier and Penot is essential. 
\par

The main tools for proving our results are presented in Section $2$. In Section $3$, we apply the tools to prove our main result. In the last section we present an illustrative application.
\section{Preliminaries}

The notation and terminology used in this paper are standard. For convenience of the reader, in this section we recall some basic facts. Let $\mathcal{F}(\Omega)$ be a family of real-valued functions defined in a set $\Omega$ arbitrary. 
\par
In sequel we need of the following definition.

\paragraph{Definition.} A mapping $T:\mathcal{F}(\Omega)\to\mathcal{F}(\Omega)$ is called strongly nonexpansive on a subset $K$ of $\mathcal{F}(\Omega)$ if 
$$
|T(u)(x)-T(v)(x)|\leq |u(x)-v(x)|,
$$
for every $x\in\Omega$ and all $u,v$ in $K$.
\par
In what follows, we will denote by 
 $L\sp{\infty}(\Omega)$ the space $L\sp{\infty}(\Omega,\Sigma,\mu)$, 
where $(\Omega,\Sigma,\mu)$ is a positive and finite measure space. In the particular case when $\mathcal{F}(\Omega)=L\sp{\infty}(\Omega)$, the inequality above must be considered in the sense almost everywhere. Clearly, every strongly nonexpansive map is nonexpansive.
\par
Our main result is as follows
\begin{theorem}\label{trm:1} Let $(\Omega,\Sigma,\mu)$ be a positive finite measure space. Then, every weakly compact convex subset of $L\sp{\infty}(\Omega,\Sigma,\mu)$ has the fixed point property for strongly nonexpansive mappings.
\end{theorem}

In the proof of the Theorem \ref{trm:1} we will use the following characterization due to Zolezzi \cite{6} of weak convergence for the space $L\sp{\infty}(\Omega,\Sigma,\mu)$.
\begin{theorem} Let $(\Omega,\Sigma,\mu)$ be a positive totally finite measure space. If an ordinary sequence $u_n\rightharpoonup u$ in $L\sp{\infty}(\Omega,\Sigma,\mu)$, then $u_n\to u$ in every $L\sp p(\Omega,\Sigma,\mu)$, $1\leq p<\infty$.
\end{theorem}  
 
Next we recall the following result from the literature on functional analysis \cite{7} which also will be used as a crucial tool to prove the Theorem \ref{trm:1}.

\begin{theorem} (Eberlein \v{S}mulian). Suppose K is weakly closed in a Banach space $X$. Then the following are equivalent:
\begin{itemize}
\item[\rm{(i)}] $K$ is weakly compact.
\item[\rm{(ii)}] $K$ is weakly sequentially compact, i.e, any sequence in $K$ has a subsequence which converges weakly.
\end{itemize}
\end{theorem}

\section{Proof}
In this section we prove the Main result. Our approach is very elementary and self-contained: In fact, beyond the theorems of Zolezzi and Eberlein-\v{S}mulian we do not use any special results other than some basic tools in the measure theory.
\par
Let $K$ be a weakly compact convex subset of $L\sp{\infty}(\Omega)$. Translating the set $K$, we may assume that $0\in K$. Now, let $T$ be a strongly nonexpansive mapping $K$ into itself. By Banach contraction principle for each $n\in\mathbb{N}$ there exists uniquely $u_n\in K$ such that

\begin{equation}\label{eqn:1}
u_n=\lambda_nT(u_n),\quad n\in\mathbb{N},
\end{equation}
where $\{\lambda_n\}$ is sequence of positive real numbers with $\lambda_n<1$ and $\lambda_n\to 1$, as $n\to\infty$. Since $K$ is weakly compact, by Eberlein-$\check{S}$mulian theorem there is $u\in K$ sucht that some subsequence $\{u_{n_j}\}$ of $\{u_n\}$ coverges weakly to $u$ in $L\sp{\infty}(\Omega)$. Now, in view of Zolezzi's theorem it follows that $u_{n_j}\to u$ in $L\sp 1 (\Omega)$, as $j\to\infty$. Therefore, up to a subsequence, we may suppose that $u_{n_j}\to u$ a.e in $\Omega$. In addition, it follows from (\ref{eqn:1}) that $T(u_{n_j})\to u$ a.e in $\Omega$. With this facts in mind, we claim now that $T(u)=u$ a.e in $\Omega$. Indeed, given $\varepsilon>0$ by Egoroff's theorem there exists $A_{\varepsilon}\subset\Omega$ such that $\mu (A_{\varepsilon})<\varepsilon$ and $u_{n_j}\to u$ uniformly on $A\sp c\sb{\varepsilon}$. Hence, 

\begin{equation}\label{eqn:2}
\int\sb{\Omega}|T(u)-u|d\mu\leq c\varepsilon+\int\sb{A\sp{c}\sb{\varepsilon}}|T(u)-u|d\mu,
\end{equation}
where $c>0$ is a constant. Applying now the Fatou's lemma in (\ref{eqn:2}), we have
\begin{equation}\label{eqn:3}
\int\sb{\Omega}|T(u)-u|d\mu\leq c\varepsilon +\lim\sb{j\to\infty}\int\sb{A\sp{c}\sb{\varepsilon}}|T(u_{n_j})-T(u)|d\mu
\end{equation}
By virtue of $T$ to be strongly nonexpansive and $u_{n_j}\to u$ uniformly on $A\sp c\sb{\varepsilon}$, (\ref{eqn:3}) yields 
\begin{equation}\label{eqn:4}
\int_\Omega |T(u)-u|d\mu\leq c\cdot\varepsilon.
\end{equation}
Since $\varepsilon$ can be arbitrarily small in (\ref{eqn:4}) we conclude that $T(u)=u$ a.e in $\Omega$. As desired. \quad $\square$

\paragraph{Remark.}  It is important to note that if $\Omega$ is a compact Hausdorff space then, a strongly nonexpansive mapping $T:C(\Omega)\to C(\Omega)$ is sequentially weakly continuous. This is due to the fact that weak convergence in $C(\Omega)$ implies in pointwise convergence in $\Omega$. In particular, arguing exactly as in the proof of Theorem \ref{trm:1}, a $C(\Omega)$-version of this result works. Thus, in general any convex and weakly compact subset of $C(\Omega)$ has the fixed point property for strongly nonexpansive self-mappings in $C(\Omega)$.

\section{Example}

\hspace{.2cm} In this section, we will use the $C(\Omega)$-version of Theorem \ref{trm:1} mentioned above to give a nontrivial example of a subset of $C[0,1]$ which is not weakly compact. Consequently, we show that the assumption of the weak compactness in the Theorem of Arino, Gautier and Penot \cite{9} is essential.

\begin{example}\label{ex:1} Let $C[0,1]$ be the space of continuous functions on $[0,1]$ endowed with its usual norm. Set $K=\{ u\in C[0,1] : u(0)=0, u(1)=1, \mbox{ and } 0\leq u(x)\leq 1\}$. Clearly, $K$ is closed, convex and bounded in $C[0,1]$. We claim now that $K$ is not weakly compact in $C[0,1]$. Suppose the contrary. Let define $T :K\to K$ by
$T(u)(x)=x\cdot u(x),\quad \forall x\in [0,1]$. Then, $T$
is well defined, i.e, $T$ maps $K$ into $K$. Now, it is straightforward to check that $T$ is a strongly nonexpansive mapping without fixed points. This contradicts the Remark above, and proves the claim. On the other hand, one can shows that $T$ is a weakly sequentially continous mapping. 
\end{example}

\paragraph{Remark.} The application above were inspired by an example of Kirk \cite{8}

\section*{Acknowledgements}
This work was supported by CAPES, Brazil.

\end{document}